\magnification=\magstep1
\documentstyle{amsppt}
\NoBlackBoxes
\def\var{\varepsilon}
\def\opartial{\overline\partial}
\def\ozeta{\overline\zeta}

\def\bC{\Bbb C}
\def\bZ{\Bbb Z}
\def\bR{\Bbb R}

\def\calP{\Cal P}
\def\calK{\Cal K}
\def\oz{\overline z}
\def\oZ{\overline Z}
\def\ds{\displaystyle}
\pageheight{9truein}
\topmatter
\title 
The Dolbeault Complex in Infinite Dimensions, II
\endtitle
\author
\bf L\'aszl\'o Lempert\footnote""{Research partially supported by an NSF grant.\hfill\hfill}
\endauthor
\endtopmatter
\TagsOnLeft
\baselineskip=12pt
\parskip=6pt
\subhead 1.\ \ Introduction\endsubhead

The goal of this series is to explore various aspects of the inhomogeneous Cauchy--Riemann, or $\opartial$, equation on infinite dimensional complex manifolds.
In the first paper in the series we argued for the importance of such an undertaking; we also gave rather complete results when the manifold in question is an infinite dimensional projective space, see [L].
In the present work we turn to the analytically more challenging problem of solving the $\opartial$ equation in an open subset of a Banach space, and offer a positive result in the space $l^1$.

Up to now not a single infinite dimensional Banach space and an open set therein have been proposed where the equation
$$
\opartial u=f\qquad (\opartial f=0)\tag1.1
$$
could be proved to be solvable under reasonably general conditions on $f$.
However, we must mention two important results.
First, Coeur\'e has constructed a continuously differentiable, closed $(0,1)$ form $f$ on a Hilbert space for which
$\opartial u=f$ has no solution on any nonempty open set, see [C,M].
Second, Raboin in [R] proved the solvability of (1.1) on the level of $(0,1)$ forms if not on open but on somewhat smaller subsets of a Hilbert space.
This result has consequences on solving (1.1) in more complicated locally convex vector spaces; in particular it implies that in a nuclear space $(1.1)$ is locally solvable for a smooth closed $(0,1)$ form $f$.
(Note though that infinite dimensional Banach spaces are never nuclear.)
For further results consult [D,L].

Why is the study of the $\opartial$ equation in infinite dimensions so much harder than in $\bC^n$?
We can discern two reasons.
One has to do with the fact that in $\bC^n$ solutions of (1.1) are constructed using integrals, and most of the time integrals with respect to volume measure; here translation invariance of the Lebesgue measure is essential.
In infinite dimensions translation invariant measures are all pathological, hence tools different from integrals must be found.

The other circumstance that complicates the study of (1.1) is that in infinite dimensions the gap between local and global is much wider than in finite dimensions.
For example, a smooth function or form on a ball in some Banach space need not be bounded on a concentric ball of somewhat smaller radius.
This latter complication we have not been able to overcome, and in our Theorem 1.1 below we had to impose a global condition, rather than the local condition of mere regularity, on the right hand side to obtain a solution of (1.1).
The global condition involved is Lipschitz continuity, whose precise meaning will be explained in Section 2, along with other basic notions of infinite dimensional analysis.

\proclaim{Theorem 1.1} Let $B(R)$ denote the open ball of radius $R\leq \infty$ centered at $0$ in the Banach space $l^1$, and suppose that $f$ is a complex valued $(0,1)$ form on $B(R)$ that is Lipschitz continuous on all balls $B(r),\ r < R$.
If $f$ is $\opartial$-closed then there is a continuously differentiable function $u$ on $B(R)$ that solves $\opartial u=f$.
If, in addition, $f$ is $m$ times continuously differentiable, $m=2,3,\ldots$, then so is $u$.
\endproclaim

At the price of a little extra effort one can prove a vector valued version of the above result.
Thus Theorem 1.1 is true for $f,\ u$ taking values in an arbitrary Fr\'echet space.

Two aspects of the hypotheses in this theorem are worth discussing.
First, Lipschitz continuity might seem too strong an assumption to the reader with experience in the finite dimensional theory.
However, in contrast with finite dimensions, in $l^1$ mere continuity of $f$ is not sufficient, see section 9.

Second, the question arises:\ why $l^1$ of all Banach spaces?
We cannot fully answer this question; and indeed similar theorems
are likely to hold for whole classes of Banach spaces (however, see
Theorem 9.1).
It is nevertheless possible to indicate how our approach depends on working in $l^1$ rather than in some other $l^p$ space, say.
It happens that we use the structure of $l^1$ in an essential way at two steps in the proof:\ in Lemma 4.1 and in Proposition 6.2.
Of these two the first one strikes us as being the more important.
It has to do with the convergence of a certain series in infinitely many variables.
The series in question is related to the geometric series
$$
\sum_{k_1=0}^\infty \sum_{k_2=0}^\infty\ldots |z_1|^{k_1} |z_2|^{k_2}\ldots |z_\nu|^{k_\nu}\ldots=\prod^\infty_{\nu=1}\sum^\infty_{k_\nu=0}|z_\nu|^{k_\nu}=\prod^\infty_{\nu=1} {1\over 1-|z_\nu|},\tag1.2
$$
if $|z_\nu| < 1$, and the latter converges only if $z=(z_\nu)\in l^1$.
The series that is of relevance to our argument dominates (1.2) termwise, and one can show that it converges if and only if $z=(z_\nu)$ is in the unit ball of $l^1$.
This is then the principal reason why we have to restrict ourselves to the space $l^1$.

Interestingly, the special role that $l^1$ plays in infinite dimensional complex analysis was also noted some ten years ago by Ryan in connection with the problem of monomial expansions of holomorphic functions, see [Ry].
That the two seemingly unrelated problems lead to the same space $l^1$ can be easily explained.
Indeed, technicalities aside, in our approach we use a torus action on $l^1$, whose holomorphic eigenfunctions are the very same monomials that occur in the series Ryan is investigating in his paper.
Also, monomial expansions themselves will appear in our proof.

From the perspective of real analysis the space $l^1$ stands out as well, but in a different sense.
According to Kurzweil the only continuously differentiable function on $l^1$ with bounded support is the zero function; in particular, there are no smooth partitions of unity, see [K].
This circumstance is regrettable since results on the $\opartial$ equation like Theorem 1.1 are often used in conjunction with partitions of unity.
Therefore it would be of great importance to explore the solvability of the $\opartial$ equation in Banach spaces that do admit smooth partitions of unity, in particular in Hilbert spaces (on this account see [DGZ]).

A few words about the proof of Theorem 1.1, which consists of two parts.
In the first, harder part we solve (1.1) on balls $B(r),\ r < R$, see Theorem 7.1.
Here we use the torus action alluded to above, and the corresponding Fourier expansions.
On the level of Fourier components (1.1) is easily solved, and the problem becomes synthesizing the componentwise solutions into an overall solution.
This is related to convergence of trigonometric series, a complicated issue in the case of infinitely many -- or even finitely, but arbitrarily many -- variables.
Indeed, no reasonable smoothness condition seems to imply pointwise convergence of Fourier series of the type considered in section 5.
We overcome this problem by using Ces\`aro--Fej\'er summation.
However, the main point is not so much convergence or summability, but
to prove a priori estimates when convergence is automatic, e.g. when
all but finitely many Fourier components vanish. We obtain the
required estimates by comparing our series with carefully selected
low dimensional solutions of (1.1) in section 6.

Once (1.1) is known to be solvable on all balls $B(r)$ with $r < R$, one constructs a solution on $B(R)$ via an approximation theorem for holomorphic functions, as in the finite dimensional case.
But there is a difference.
In $\bC^N$ one can use the fact that given a holomorphic function $\varphi$ on $\{\|z\| < R\},\ r < R$, and $\var > 0$, there is an entire function $\psi$ that approximates $\varphi$ within $\var$ on $\{\|z\| < r\}$.
A similar uniform approximation theorem is not known to hold in any infinite dimensional Banach space; yet the day can be saved by introducing a norm different from sup norm, in which approximation is possible.
The last step of the proof then is an application of this approximation theorem to solve (1.1) on $B(R)$, see section 8.

\subhead 2.\ \ Background\endsubhead

Although in [L] we introduce the formalism of the $\opartial$ equation in complete generality, for the reader's convenience we quickly review those aspects that will be of relevance to this paper.
Thus suppose $V$ is a complex Banach space, $\Omega\subset V$ is open, and $u\colon\Omega\to \bC$ is a function.
For $z\in\Omega,\ \xi\in V$ we put
$$
du(z;\xi)=\lim(u(z+t\xi)-u(z))/t\quad ,\ \bR\ni t\to 0,\tag2.1
$$
and say that $u$ in continuously differentiable, or $u\in C^1(\Omega)$, if the limit in (2.1) exists for all $z\in\Omega,\ \xi\in V$, and $du\colon\Omega\times V\to \bC$ is continuous.
In general, we inductively define $u\in C^m(\Omega)$ if $du\in C^{m-1}(\Omega\times V),\ m=1,2,\ldots,\ C^0(\Omega)=C(\Omega)$.
In addition, for any set $A\subset\Omega$ we introduce the quantities
$$
\aligned
|u|_{0,A}&=\sup_{z\in A}|u(z)|\\
|u|_{1,A}&=\sup_{z\not= w\in A}|u(z)-u(w)|/\|z-w\|.
\endaligned\tag2.2
$$
If $|u|_{1,A} < \infty$ we say that $u$ is Lipschitz continuous on $A$.
When $A=\Omega$ we simply write $|u|_0,\ |u|_1$ for $|u|_{0,A},\ |u|_{1,A}$.

A 1--form on $\Omega$ is a function $f\colon\Omega\times V\to \bC$ that is $\bR$--linear in the $V$--variable; if it is also $\bC$--linear resp.~$\bC$--antilinear, we say $f$ is a $(1,0)$ resp.~$(0,1)$ form.
If $u\in C^1(\Omega)$ then $du$ is a 1--form, which can be uniquely decomposed into the sum $\partial u+\opartial u$ of a $(1,0)$ and a $(0,1)$ form; this decomposition defines the $\partial$ and $\opartial$ operators on $C^1$ functions.
We say that a 1--form $f$ is of class $C^m$ if $f$ is $C^m$ as a function on $\Omega\times V$, and we denote the space of $(0,1)$--forms of class $C^m$ on $\Omega$ by $C^m_{0,1}(\Omega)$.
With $B$ denoting the unit ball in $V$, $f$ a 1--form on $\Omega$, and $A\subset \Omega$ we put
$$
\aligned
|f|_{0,A}&=\sup_{(z,\xi)\in A\times B} |f(z;\xi)|,\\
|f|_{1,A}&=\sup_{z\not= w\in A,\xi\in B} |f(z;\xi)-f(w;\xi)|/\|z-w\|.\\
\endaligned\tag2.3
$$
We say that $f$ is Lipschitz continuous on $A$ if $|f|_{0,A} < \infty$ and $|f|_{1,A} < \infty$, and when $A=\Omega$ for $|f|_{0,A},\ |f|_{1,A}$ we simply write $|f|_0,\ |f|_1$.

More generally, we can define the notions above for vector valued functions or forms.
We shall only need this when $u,f$ take values in a Banach space, in fact in the same space $V$.
In this case (2.1) still can be used to define differential and smoothness classes, and (2.2), (2.3), with absolute values replaced by norms, will define $|u|_{m,A},\ |f|_{m,A},\ m=0,1$.
We shall also need to know how the quantities $|f|_0, |f|_1$ transform under $C^1$ mappings $\sigma\colon\tilde\Omega\to \Omega$.
A simple computation gives
$$
\aligned
|\sigma^* f|_0 &\leq |\sigma|_1 f|_0,\\
|\sigma^* f|_1 &\leq |\sigma|^2_1 |f|_1+|d\sigma|_1 |f|_0.
\endaligned\tag2.4
$$
\comment
More generally, we shall have to deal with vector valued functions and forms.
Let $W$ we a locally convex (topological vector) space over $\bC$ whose topology is defined by a family $\calP$ of seminorms $p\colon W\to\bR$.
Without further mentioning, we shall always assume the $W$ is sequentially complete.
For a function $u\colon\Omega\to W$ its differential $du$ is again defined by (2.1), and the smoothness classes $C^m(\Omega,W)$ are defined inductively as before.
Similarly, we have the classes $C^m_{0,1}(\Omega,W)$ of $m$ times continuously differentiable $W$ valued $(0,1)$ forms $f$, and the operator $\opartial\colon C^m(\Omega,W)\to C^{m,1}_{0,1}(\Omega,W),\ m=1,2,\ldots$.
If $p\in\calP$ and $A\subset\Omega$, the quantities $|u|^{(p)}_{0,A},\ |u|^{(p)}_{1,A},\ |f|^{(p)}_{0,A}$, etc.~are introduced by an obvious modification of (2.2), (2.3); for example
$\ds{u^{(p)}_{0,A}=\sup\limits_{z\in A} p(u(z))}$.
A function $u$ resp.~form $f$ is said to be bounded or Lipschitz continuous on $A$ if for all $p\in\calP\ |u|^{(p)}_{0,A} < \infty$ resp.~$|f|^{(p)}_{0,A} < \infty$ or $|u|_{1,A}^{(p)} < \infty$ resp.~$|f|^{(p)}_{1,A} < \infty$.
With these notions we shall be able to prove Theorems 1.1, 1.2 for $W$ valued forms, with the proviso that in case of Theorem 1.2 we shall have to assume that the family $\calP$ of seminorms is countable, i.e.~$W$ is a Fr\'echet space.
\endcomment

The most useful topology to endow the space of 
functions or forms with is the compact open topology.
In this topology $u_j\to u$ resp.~$f_j\to f$ if uniform convergence holds on all compact subsets of $\Omega$ resp.~of $\Omega\times V$.
Thus $u_j\to u$ is equivalent to $|u_j-u|_{0,K}\to 0$ for all compact $K\subset\Omega$ (however, $f_j\to f$ is weaker than $|f_j-f|_{0,K}\to 0$).
The following simple proposition will be very useful:

\proclaim{Proposition 2.1} If $u_j\in C(\Omega)$ and $u_j\to u$ uniformly on compact subsets of $\Omega$ then $u\in C(\Omega)$.
\endproclaim

\demo{Proof} Clearly the restriction of $u$ to any compact $K\subset\Omega$ is continuous.
If $z_n\in\Omega$ converges to $z\in\Omega$ then $K=\{z,z_1,z_2,\ldots\}$ is compact, therefore $u(z)=\lim u(z_n)$:\ \ thus $u$ is continuous.
\enddemo

We have already introduced the equation $\opartial u=f$ for $u\in C^1(\Omega),\ f\in C_{0,1}(\Omega)$, next we will extend its meaning.
Suppose $u\in C(\Omega),\ f\in C_{0,1}(\Omega)$.
We say that $\opartial u=f$ (in the weak sense) if for any finite dimensional subspace $F\subset V\ \opartial (u|_F)=f|_F$ holds in the sense of distribution theory.
(Considering only one dimensional subspaces $F$ would give the same notion.)
When $V$ itself is finite dimensional, it is easy to verify that $\opartial u=f$ in the above sense is equivalent to the equation in the sense of distribution theory.
A straightforward consequence of the definition and of Proposition 2.1 is

\proclaim{Proposition 2.2} If $u_j\in C(\Omega),\ f_j\in C_{0,1}(\Omega),\ \opartial u_j=f_j$, and $u_j\to u$ resp.~$f_j\to f$ uniformly on compact subsets of $\Omega$ resp.~$\Omega\times V$ then $\opartial u=f$.
\endproclaim

We shall also need the following facts from ``elliptic regularity theory'':

\proclaim{Proposition 2.3} If $u\in C(\Omega),\ f\in C_{0,1}^m (\Omega),\ m=1,2,\ldots$, (resp.~$f$ is Lipschitz continuous) and $\opartial u=f$ (weakly) then $u\in C^m(\Omega)$ (resp.~$u\in C^1(\Omega)$) and $\opartial u=f$ holds according to the original definition of $\opartial$.
\endproclaim

\proclaim{Proposition 2.4} If $u_j\in C(\Omega)$ and $\opartial u_j\in C_{0,1}(\Omega)$ are uniformly bounded on compact subsets of $\Omega$ resp.~$\Omega\times V$ then the $u_j$ are locally equicontinuous on $\Omega$ (i.e.~$\Omega$ can be covered by open subsets on which the $u_j$ are equicontinuous).
\endproclaim


For the proof of the first case of Proposition 2.3 see [L, Proposition 9.3].
Proposition 2.4 follows by first noticing that the assumptions imply that $u_j, \opartial u_j$ are locally uniformly bounded; second, by restricting to one dimensional slices, and applying Pompeiu's representation formula (also known as Cauchy formula, see [H\"o, Theorem 1.2.1]) on these slices.
The same representation formula gives the second case of Proposition 2.3.

In the same spirit we shall say that $f\in C_{0,1}(\Omega)$ is $\opartial$ closed, and write $\opartial f=0$, if for all finite dimensional subspaces $F\subset V\ \opartial (f|_F)=0$ in the sense of distribution theory.
Clearly, for the equation $\opartial u=f$ to be solvable by $u\in C(\Omega)$ it is necessary that $\opartial f=0$.
A continuous function $u$ is holomorphic if $\opartial u=0$; by Proposition 2.3 this implies $u$ is $C^\infty$. The set of holomorphic functions on
$\Omega$ is denoted $\Cal O(\Omega)$.

In this paper our main concern will be the Banach space $l^1$, and so, unless indicated otherwise, $\|\ \ \|$ will denote the $l^1$ norm on $l^1$ or on $\bC^N$:\ \ if $z=(z_\nu)\in l^1$ or $\bC^N$, $\|z\|=\Sigma |z_\nu|$.
$B(R)$ resp.~$B_N(R)$ will denote the ball $\|z\| < R$ in $l^1$
resp.~$\bC^N$.
Unless indicated otherwise, $R < \infty$.

\subhead 3.\ \ The equation $\opartial u=f$ in $\bC^N$\endsubhead

An obvious approach to solving the equation
$$
\opartial u=f\tag3.1
$$
in infinite dimensions, which has been attempted over and over, is to solve (3.1) in $\bC^N$, obtain pointwise (or preferably uniform) estimates for the solution $u$, and see what happens when $N\to\infty$.
If the estimates are independent of $N$, it is reasonable to expect that from the finite dimensional solutions one can construct, by some limiting procedure, a solution of the infinite dimensional equation.
However, such a direct approach fails, because the finite dimensional estimates for the solution of (3.1) tend to blow up as $N\to\infty$, and exponentially
at that.
Below we shall derive one such estimate; it is not the sharpest, but the one that seems easiest to present.
In spite of its shortcomings, it will serve as our starting point.

Let $T^N=\bR^N/\bZ^N$ denote the $N$--dimensional torus.
It acts on $\bC^N$ by
$$
\rho_t(z_1,\ldots,z_N)=(e^{2\pi i t_1}z_1,\ldots, e^{2\pi i t_N} z_N)\quad,
\quad
(t_\nu)=t\in T^N.
$$
{\it In this section $\|\quad\|$ will denote any norm on $\bC^N$ that is invariant under the action $\rho$, and} $B_N=\{z\in\bC^N\colon\|z\| < 1\}$.
In particular, $\|\quad\|$ can be the $l^1$ norm.

Let $\|\quad\|^*$ denote the dual norm on $\bC^N$ and $\alpha$ a positive number such that
$$
(\Sigma_\nu |w_\nu|^2)^{1/2}\leq \alpha \|w\|^*\qquad ,\ \ w=(w_\nu).
$$
Let diam $\Omega$, Vol $\Omega$ denote the {\it Euclidean} diameter and volume of a set $\Omega\subset\Bbb C^N$, and $r(z)$ the distance of a point
$z\in\Omega$ to $\partial\Omega$, measured in the norm $\|\quad\|$.

\proclaim{Lemma 3.1} Suppose that $f$ is a complex valued, bounded and
measurable  $(0,1)$ form on some bounded pseudoconvex open set $\Omega\subset\bC^N$.
If $\opartial f=0$ then there is a $u\in C(\Omega)$ such that $\opartial u=f$ and for $z\in\Omega$
$$
|u(z)|\leq \{2\alpha\ \text{\rm diam }\Omega(\text{\rm Vol }\Omega/\text{\rm Vol }B_N)^{1/2} r(z)^{-N}+2r(z)\}|f|_0.\tag3.2
$$
\endproclaim

In the lemma $\opartial$ refers to the Cauchy--Riemann operator extended to act on distributions.

\demo{Proof} Write $f=\Sigma f_\nu dz_\nu$; then
$$
\sup_\Omega (\Sigma |f_\nu (z)|^2)^{1/2}\leq \alpha\sup_\Omega\|(f_\nu (z))\|^*=\alpha|f|_0.
$$
Hence by a theorem of H\"ormander (cf.~[H\"o, Theorem 4.4.2]) there is a square 
integrable solution $u\colon \Omega\to \bC$ of the equation $\opartial u=f$ such that
$$
\int_\Omega |u|^2\leq \alpha^2 (1+\text{ diam}^2\Omega)^2\text{ Vol }\Omega |f|_0^2.
$$
A well known trick turns this into the existence of a solution $u$ such that
$$
\int_\Omega |u|^2\leq (2\alpha)^2\text{ diam}^2\Omega\text{ Vol }\Omega |f|_0^2.\tag3.3
$$
Indeed, when diam $\Omega=1$, this is what we had before; a general $\Omega$ can be scaled down to one with diameter 1.
Since under scaling $\int_\Omega |u|^2$ transforms as Vol $\Omega$
while $|f|_0$ diam $\Omega$ does not change, (3.3) is obtained generally.

In fact the solution $u$ is even continuous.
Indeed, by [GL,H] there are continuous local solutions $u'$ of $\opartial u'=f$; since $u-u'$ is holomorphic, $u$ itself must be continuous. 

With $z\in\Omega,\ r=r(z)$, and $d\lambda=dt_1\ldots dt_N$ denoting Lebesgue measure on $T^N$ --- thus $\lambda(T^N)=1$ --- define a function $v$ on $B_N$ by
$$
v(\zeta)=\int_{T^N} u(z+r\rho_t(\zeta))d\lambda(t).
$$
We will estimate the Lipschitz norm of $v$, first by assuming
$v\in C^1(B_N)$.
Then $|\opartial v|_0\leq r|f|_0$.
Also $v$ is invariant under the action $\rho$, which implies that
$dv(z,\xi)=0$ if $z\in B_N\cap\bR^N$, $i\xi\in\bR^N$. It follows that
$dv|_{\bR^N}=2\opartial v|_{\bR^N}$ and so for $\zeta,\zeta'\in B_N\cap \bR^N$
$$
|v(\zeta)-v(\zeta')|\leq 2|\opartial v|_0\|\zeta-\zeta'\|\leq
2r|f|_0 \|\zeta-\zeta'\|. \tag3.4
$$
In fact, by invariance (3.4) holds for all $\zeta,\zeta'\in B_N$; and a simple
approximation argument allows us to extend (3.4) to the general
case, without the assumption $v\in C^1(B_N)$. Furthermore
$$
\int_{B_N} |v|^2\leq r^{-2N}\int_\Omega |u|^2.\tag3.5
$$
Choose a point $\zeta\in B_N$ where $|v|^2$ is not greater than its average on $B_N$; by (3.3), (3.5)
$$
|v(\zeta)|\leq 2\alpha\text{ diam }\Omega (\text{Vol }\Omega/\text{Vol }B_N)^{1/2} r^{-N} |f|_0.
$$
Finally apply (3.4) with $\zeta'=0$ to obtain the required estimate (3.2) for $u(z)=v(0)$.
\enddemo

An easy computation gives

\proclaim{Corollary 3.2} In the special case when $\Omega=\{z\in\bC^N\colon \|z\| < R\},\ 0 < R < \infty$, the solution $u$ satisfies
$$
|u(z)|\leq 2(R+\alpha\text{\rm diam }\Omega)
\left({R\over R-\|z\|}\right)^N |f|_0.
$$
When $\|\quad\|$ is the $l^1$ norm, {\rm diam} $\Omega=2R$, and $\alpha$ can be chosen $\sqrt{N}$.
\endproclaim

\subhead 4.\ \ A series in infinitely many variables\endsubhead

Much finite dimensional analysis depends on the convergence of the geometric series.
Similarly, we shall base our analysis in $l^1$ on a series in infinitely many variables.

If $z=(z_\nu)$ is a finite or infinite sequence of numbers, we shall denote by \#$z$ the number of nonzero entries $z_\nu$.
Thus $0\leq \#z\leq \infty$.
Further, we shall write $z\geq 0$ to indicate that all $z_\nu$ are real and nonnegative.
If $z=(z_\nu)$ is in the unit ball $B(1)$ of $l^1$ and $q\in \Bbb C$, put
$$
\Delta(q,z)=\sum_{k=(k_\nu)\geq 0} {(\sum_\nu k_\nu)^{\sum\limits_\nu k_\nu}\over \prod_\nu k_\nu^{k_\nu}}\ |q|^{\# k}|z^k|.\tag4.1
$$
The first summation is extended over integer multiindices $k=(k_\nu)^\infty_{\nu=1}\geq 0$ such that $\#k < \infty$.
{\it In this section $k$ will always denote such a multiindex}.
In (4.1) $z^k$ stands for $\prod_\nu z_\nu^{k_\nu}$, and we use the convention $0^0=1$.
Below we shall occasionally write the coefficients above in the more compact form $\|k\|^{\|k\|}/k^k$.

\proclaim{Lemma 4.1} The series (4.1) is uniformly convergent on
compact subsets of $\Bbb C\times B(1)$ and $\Delta$ is continuous on
$\Bbb C\times B(1)$.
\endproclaim

Convergence of (4.1) with $q=1$ is the issue, albeit implicitly, in Ryan's
paper [Ry] mentioned in the Introduction.
His estimates amount to proving that $\Delta(1,z)$ is bounded when $\|z\|\leq 1/(e+\varepsilon)$ for any $\var > 0$.

\demo{Proof} We shall use the following inequalities:
$$
\align
& s!\leq s^s\leq e^s s!\qquad ,\ \ s=0,1,2,\ldots;\tag4.2\\
& (e \eta)^s s!\leq s^s\qquad\quad ,\ \ s\geq s_\eta,\tag4.3
\endalign
$$
for any given $\eta < 1$, provided $s_\eta$ is sufficiently large.
(4.3) follows from Stirling's formula, (4.2) is easily proved by induction.

With $Q\geq 1$, $0<\eta < 1$, and $n=0,1,\ldots$ consider the set
$$
H(Q,\eta,n)=
\{(q,z)\in \overline{B_1(Q)}\times\overline{B(\eta^3)}\colon \sum_{\nu > n} |z_\nu|\leq (\eta-\eta^2)/(e|q|)\}.
\tag4.4
$$
The sets int $H(Q,\eta,n)$ cover $\Bbb C\times B(1)$.
First we will estimate $\Delta$ uniformly on $H(Q,\eta,n)$.

If $(q,z)\in H(Q,\eta,n)$, define $X(z)=X=(X_\nu)$ by
$$
X_\nu=\cases |z_\nu|/\eta&\text{if $\nu\leq n$}\\
e|qz_\nu|&\text{if $\nu > n$}\endcases\tag4.5
$$
and check that $\|X\|\leq \eta$.
Next choose a multiindex $k$ in (4.1), and with $s_\eta$ as in (4.3) denote by $F$ the set of those $\nu\in\{1,\ldots,n\}$ that satisfy $k_\nu > s_\eta$.
Using (4.2), (4.3), and $\#k\leq n+\sum\limits_{\nu > n}k_\nu$ we can estimate
$$
\align
{\left(\sum_\nu k_\nu\right)^{\sum\limits_\nu k_\nu}\over
\prod_\nu k_\nu^{k_\nu}} &|q|^{\# k}|z^k|\leq { e^{\sum\limits_\nu k_\nu}
(\sum_\nu k_\nu)! Q^{n+\sum\limits_{\nu > n} k_\nu}\over \prod\limits_{\nu\not\in F}k_\nu!\prod\limits_{\nu\in F} (e\eta)^{k_\nu} k_\nu!}|z^k|\\
&\leq Q^n {\left(\sum_\nu k_\nu\right)!\over \prod_\nu k_\nu!} e^{\sum\limits_{\nu\not\in F} k_\nu} \eta^{-\sum\limits_{\nu\leq n}k_\nu} Q^{\sum\limits_{\nu > n} k_\nu} |z^k|\\
&\leq Q^n e^{n s_\eta} {\left(\sum_\nu k_\nu\right)!\over \prod_\nu k_\nu!}\eta^{-\sum\limits_{\nu\leq n} k_\nu}(eQ)^{\sum\limits_{\nu > n} k_\nu}|z^k|\\
&\leq Q^n e^{n s_\eta} {\left(\sum_\nu k_\nu\right)!\over
\prod_\nu k_\nu!} X^k.
\endalign
$$
Thus
$$
\aligned
\Delta (q,z)&\leq Q^n e^{n s_\eta}\sum^\infty_{j=0}\ \ \sum_{k\geq 0,\Sigma k_\nu=j} {j!\over \prod k_\nu!}\ X^k\\
&\leq Q^n e^{n s_\eta} \sum^\infty_{j=0} \|X\|^j \leq Q^n e^{n s_\eta}/(1-\eta)
\endaligned\tag4.6
$$
shows that $\Delta$ is bounded on $H(Q,\eta,n)$.

The lemma now follows from
Proposition 4.2 below and Proposition 2.1.
\enddemo

\proclaim{Proposition 4.2}
If $\Omega$ is an open set in any Banach space $V$, and $h_j\in\Cal O(\Omega)$
is a sequence such that $\sum_1^\infty|h_j|$ is bounded on $\Omega$, then this series is uniformly convergent on compacts.
\endproclaim

\demo{Proof} With arbitrary complex numbers $\var_j$ of modulus one the partial sums $\sum_1^N \var_j h_j$ are uniformly bounded on $\Omega$, whence one--variable Cauchy estimates imply that $\sum_1^N \var_j dh_j(z;\xi)$ are locally uniformly bounded on $\Omega\times V$.
The bounds being independent of $\var_j, N$, it follows that $\sum^N_1|dh_j(z;\xi)|$ are bounded, locally uniformly in $z,\xi$, and uniformly in $N$.
Thus the partial sums $\sum_1^N|h_j|$ are locally uniformly equicontinuous, which easily implies that $\sum^\infty_1|h_j|$ is indeed uniformly convergent on compacts.
\enddemo

\proclaim{Corollary 4.3} For any $\theta < 1$ there is a constant $c > 1$ such that if $\|z\| < \theta$ and $\#z < \infty$ then
$$
\Delta (q,z)\leq \max (1,|q|^{\#z}) e^{c\# z}.\tag4.7
$$
\endproclaim

\demo{Proof} Both sides are 1 when $\#z=0$.
Otherwise we can assume the first $n=\#z$ coordinates of $z$ are nonzero, and 
also $|q|\geq 1$. Thus $(q,z)\in H(|q|,\theta^{1/3},n)$, which
implies (4.6) and so (4.7) if $c$ is sufficiently large.
\enddemo

Later we shall have to deal with series of form (4.1) but with $z=(z_\nu)_1^N$ a finite sequence.
We shall put $\Delta(q,z)=\Delta(q,z_1,\ldots,z_N,0,0,\ldots)$; then all of the above, in particular (4.7) with $c$ independent of $N$, continues to hold.

Lemma 4.1 will be used repeatedly in this paper.
For the moment we want to point out some consequences pertaining to monomial expansion of holomorphic functions.
Let
$T=\prod_1^\infty (\bR/\bZ)$
denote the infinite dimensional torus, a compact topological group, and denote by $\lambda$ the Haar measure on $T,\ \lambda(T)=1$.
Monomial expansions arise from the continuous action of $T$ on $l^1$ 
$$
\rho_t(z_1, z_2,\ldots)=(e^{2\pi it_1} z_1,\ e^{2\pi i t_2} z_2,\ldots),\quad z=(z_\nu)\in l^1,\ t=(t_\nu)\in T.
$$
Given $0 < R \leq \infty$ and  $h\in\Cal O(B(R))$, its monomial expansion is
$$
h\sim\sum_k h_k\quad ,\ h_k=\int_T e^{-2\pi i kt}\rho_t^* h d\lambda (t).\tag4.8
$$
The terms $h_k\in\Cal O (B(R))$ transform under $\rho$ as
$\rho_t^* h_k=e^{2\pi i kt} h_k$.
Upon restricting to finite dimensional coordinate planes $P$ one finds that $h_k$ is indeed a monomial of form $a_k z^k,\ a_k\in \bC$, and the monomial expansion of $h$ becomes
$$
h(z)\sim\sum a_k z^k.\tag4.9
$$
Clearly the restriction of this series to coordinate planes $P$ as above is just the Taylor series of $h|_P$, hence $\Sigma a_k z^k=h(z)$ at least when $z\in B(R),\ \#z < \infty$.

\proclaim{Theorem 4.4} (a) If $h$ is a bounded holomorphic function on $B(R),\ R < \infty$, then
$$
a=\sup_k |a_k| R^{\|k\|} k^k/\|k\|^{\|k\|} < \infty,\tag4.10
$$
and the series (4.9) converges to $h$.

(b) Conversely, if a sequence $\{a_k\}$ satisfies (4.10) then the series $\Sigma a_k z^k$ converges, uniformly absolutely on compact subsets of $B(R)$, to a
$g\in\Cal O(B(R))$ whose monomial expansion is $\Sigma a_k z^k$.
Furthermore, for any compact $K\subset B(R)$ there is a constant $C_K$ such that
$$
\max _K |g|\leq C_K a.\tag4.11
$$
\endproclaim

\demo{Proof} (b) The series $\Sigma |a_k| |z^k|$ is termwise dominated by the series $a\Delta (1,z/R)$, so that by Lemma 4.1 $\Sigma a_k z^k$ is indeed uniformly absolutely convergent on compact subsets of $B(R)$; its sum $g$ is holomorphic by Proposition 2.2, and (4.11) holds with $C_K=$ max$_K |\Delta (1,z/R)| < \infty$.
Computing the integral in (4.8) with $h=g$, term by term, gives $h_k=a_k z^k$, as claimed.

(a) Let $0 < r < R$ and $M=\sup |h|$.
(4.8) implies $|h_k(z)|\leq M$ for $z\in B(R)$.
Putting $z=rk/\|k\|$ if $k\not= 0$, we obtain
$$
|a_k| r^{\|k\|} k^k/\|k\|^{\|k\|}\leq M;
$$
thus letting $r\to R$ (4.10) follows.
Part (b) implies the series (4.9) converges to a holomorphic function on $B(R)$, which must agree with $h$ since the two agree where $\#z < \infty$.
\enddemo

If now $0 < R\leq \infty$, and $h$ is any holomorphic function on $B(R)$, following Ryan we observe that with an arbitrary sequence $\sigma=(\sigma_1,\sigma_2,\ldots)$ such that
$$
0\leq \sigma_\nu < R\quad ,\qquad \lim_\nu \sigma_\nu=0,\tag4.12
$$
the function $h(\sigma_1 z_1,\sigma_2 z_2,\ldots)$ is bounded and holomorphic for $z\in B(1)$.
A reasoning as in [Ry], coupled with Theorem 4.4 then gives the following

\proclaim{Theorem 4.5} (a) The monomial expansion of any
$h\in\Cal O(B(R)),\ R\leq\infty$, converges to $h$, uniformly absolutely on compact subsets of $B(R)$.
The monomial coefficients $a_k$ satisfy 
$$
\text{\rm{lim}}_k |a_k| \sigma^k k^k/\|k\|^{\|k\|}=0\tag4.13
$$
for any $\sigma=(\sigma_\nu)$ as in (4.12).

(b) Conversely, if (4.13) holds for any sequence $\sigma=(\sigma_\nu)$ as in (4.12) then $\Sigma a_k z^k$ is the monomial expansion of a function holomorphic on $B(R)$.
\endproclaim

Theorem 4.5 is proved in [Ry] in the case when $R=\infty$.

\subhead 5.\ \ The key lemma\endsubhead

We again reserve $\|\quad\|$ to denote the $l^1$ norm (in this section on $\bC^N$), and consider a complex valued, bounded and measurable
$(0,1)$ form $f$ on a ball $B_N(R)\subset \bC^N,\ \opartial f=0$.
Given $0 < r < R$, we shall construct a canonical solution $u\in C(B_N(r))$ of the equation $\opartial u=f$ that is comparable to the ``smallest'' solutions in the sense explained below.
The construction will be based on Fourier series induced by the torus action $\rho$ considered in section 3.

Thus let
$$
f\sim\sum_{k\in \bZ^N}f_k\qquad ,\qquad f_k=\int_{T^N} e^{-2\pi ikt}\rho_t^* f d\lambda (t)\tag5.1
$$
be the Fourier expansion of $f$.
The series in (5.1) converges to $f$ in the sense of distribution theory; the terms $f_k$ are $\opartial$--closed, bounded and measurable $(0,1)$ forms on $B_N(R)$ that satisfy
$$
\rho_t^* f_k=e^{2\pi ikt} f_k.\tag5.2
$$
For each $k\in\bZ^N$ we can solve the equation $\opartial u_k=f_k$ with $u_k\in C(B_N(R))$, for example by Lemma 3.1.
At the price of replacing $u_k$ by $\int e^{-2\pi ikt}\rho_t^* u_k d\lambda(t)$, we can assume that $u_k$ transforms as $f_k$:
$$
\rho_t^* u_k=e^{2\pi ikt} u_k.\tag5.3
$$
This determines $u_k$ up to a holomorphic term that transforms as $u_k$
itself in (5.3).
When $k\geq 0$, this means that this holomorphic term is a constant multiple of the monomial $z^k$; for other values of $k$ such a holomorphic term vanishes, and so $u_k$ is uniquely determined.
To determine $u_k$ unambiguously for all $k$, put for $k\geq 0$
$$
Z(k)=\cases 0&\text{if $k=0$}\\
rk/\|k\|&\text{otherwise},\endcases
$$
and require that $u_k(Z(k))=0$.
Note that $Z(k)\in \partial B_N(r)$ unless $k=0$.

When $f$ is not only bounded but sufficiently many times (depending on $N$!) differentiable, $\Sigma u_k$ can be shown to converge to a solution $u$ of the equation $\opartial u=f$; yet this observation will not be of help as $N\to\infty$.
Instead, we shall show that the Ces\`aro--Fej\'er summation process alway sums $\Sigma u_k$ to a solution $u$, and this $u$ is comparable to the ``smallest'' solutions not only on $B_N(r)$, but even on $B_N(r)$ intersected with arbitrary dimensional coordinate planes.

We have to recall certain facts about Fourier series.
We start with the Fourier expansion of the iterated Fej\'er kernel
$$
\prod^N_{\nu=1} {\sin^2\pi j t_\nu\over j\sin^2 \pi t_\nu}=\prod^N_{\nu=1} \sum^{j-1}_{h=1-j} (1-{|h|\over j}) e^{2\pi iht_\nu}=\sum_{\max_\nu|k_\nu| < j} a_k^j e^{2\pi ikt},
$$
$j=1,2,\ldots$.
Given any series $\sum\limits_{k\in\bZ^N} v_k$ its iterated Ces\`aro means are
$$
v^j=\sum_{\max_\nu|k_\nu| < j} a_k^j v_k
$$
with the coefficients defined above.
In particular, if $v$ is a bounded measurable function on $B_N(R)$, with Fourier series $\Sigma v_k$, its iterated Ces\`aro--Fej\'er means are
$$
v^j=\sum_k a_k^j v_k=\int_{T^N} \prod^N_{\nu=1} {\sin^2\pi jt_\nu\over j\sin^2\pi t_\nu} \rho_t^* v d\lambda(t).
$$
The circumstance that the iterated Fej\'er kernels form a nonnegative approximation of the $\delta$ distribution implies that $v^j$ converge to $v$ as distributions, and when $v$ is uniformly continuous, the convergence is even uniform (Fej\'er's theorem).
Further, $|v^j|_m\leq |v|_m$ for $m=0,1$.
All of the above also applies to forms rather than functions.

Next, given $n\leq N$, embed $\bC^n$ into $\bC^N$ by $\bC^n \ni \zeta\mapsto (\zeta,0,\ldots,0)\in \bC^N$.
Thus $B_n(R)=B_N(R)\cap \bC^n$.
Given a $\opartial$--closed, bounded and measurable $(0,1)$ form $f$ on $B_N(R)$, and $r < R$, construct the series $\Sigma u_k$ as described above.

\proclaim{Lemma 5.1} Suppose that for some $n\leq N$ and nonnegative numbers
$A$ and $Q$ there is a $U\in C(\overline{B_n(r)})$ that solves $\opartial U=f|_{B_n(r)}$, and satisfies
$$
|U (z)|\leq A Q^{\# z},\qquad z\in\overline{B_n(r)}.\tag5.4
$$
Then the iterated Ces\`aro means of the series $\Sigma u_k|_{B_n(r)}$ locally uniformly converge to a function $u'\in C(B_n(r))$ that solves $\opartial u'=f|_{B_n(r)}$ and satisfies
$$
|u' (z)-U(z)|\leq A \Delta (Q, z/r)\quad ,\ z\in B_n(r).\tag5.5
$$
\endproclaim

\proclaim{Corollary 5.2} For any complex valued, $\opartial$--closed, bounded 
and measurable $(0,1)$ form $f$ on $B_N(R)$ the iterated Ces\`aro means of the series $\Sigma u_k$ locally uniformly converge to a solution $u\in C (B_N(r))$ of the equation $\opartial u=f|_{B_N(r)}$.
\endproclaim

\demo{Proof} By Lemma 3.1 there is a $U\in C(B_N(R))$ that solves $\opartial U=f$.
Thus the hypotheses of Lemma 5.1 are satisfied with $n=N$,
$A=\sup_{B_N(r)} |U|$, and $Q = 1$, whence the Corollary follows.
\enddemo

We shall call this $u$ the canonical solution based on the nodes $Z(k)$.
Note that $u'$ in Lemma 5.1 is nothing but $u|_{B_n(r)}$.

\demo{Proof of Lemma 5.1} Using the $T^n$ action on $\bC^n$, expand $U$ into a Fourier series.
Thinking of $\bZ^n$ as $\bZ^n\oplus (0)\subset \bZ^N$, this expansion can be written
$$
U\sim \sum_{k\in \bZ^n\oplus (0)} U_k\quad ,\ 
U_k=\int_{T^N} e^{-2\pi ikt}\rho_t^* U d\lambda (t).\tag5.6
$$

Again, $\rho_t^* U_k=e^{2\pi ikt} U_k$, and (5.4) implies
$$
|U_k (z)|\leq A Q^{\# z}.\tag5.7
$$

Since $\opartial U_k=f_k|_{B_n(r)}=\opartial u_k|_{B_n(r)}$ when $k\in \bZ^n\oplus (0)$, it follows that $U_k-u_k$ is holomorphic on $B_n(r)$.
As before, this implies $U_k-u_k|_{B_n(r)}=0$ unless $k\geq 0$; when $k\geq 0$ this difference must be a constant multiple of the monomial $z^k$.
Remembering $u_k (Z(k))=0$, we conclude $U_k(z)-u_k(z)=U_k (Z(k)) Z(k)^{-k} z^k$, $z\in B_n(r)$.
Hence by (5.7)
$$
\sum_{k\in \bZ^n\oplus (0)}|U_k(z)-u_k(z)|\leq A \sum_{k\geq 0} {(\sum_\nu k_\nu)^{\sum_{\nu} k_\nu}\over \prod_\nu k_\nu^{k_\nu}}\ Q^{\# k}|(z/r)^k|=
A\Delta (Q,z/r).\tag5.8
$$

Notice also that (5.3) implies $u_k|_{B_n(R)}=0$ if $k\in \bZ^N\backslash (\bZ^n\oplus (0))$.
Thus, by Lemma 4.1 $\Sigma U_k$ and $\Sigma u_k$ are locally uniformly equiconvergent on $B_n(r)$.
Since the former series, as the Fourier series of a continuous function, is uniformly Ces\`aro--Fej\'er summable to $U$, it follows that on $B_n(r)$ the Ces\`aro means of $\Sigma u_k$ converge, locally uniformly to a function $u'\in C(B_n(r))$.
As $U-u'$ is the locally uniform limit of holomorphic functions, we have $\opartial u'=\opartial U|_{B_n(r)}=f|_{B_n(r)}$; and (5.8) implies (5.5).
The proof is complete.
\enddemo

\proclaim{Corollary 5.3} Suppose $f_1,\ f_2$ are bounded, measurable
$\opartial$--closed $(0,1)$ forms on $B_N(R)$, and let $u_j,\ j=1,2$
denote the canonical solution of $\opartial u_j=f_j$ based on the nodes
$Z(k)$. If $f_1|_{B_n(R)}=f_2|_{B_n(R)}$ then $u_1|_{B_n(r)}=u_2|_{B_n(r)}$.
\endproclaim

\demo{Proof} If $f=f_1-f_2$, the canonical solution of $\opartial u=f$
is $u=u_1-u_2$. Since the hypotheses of Lemma 5.1 are satisfied with
$U=0$, $A=Q=0$, it follows that $u|_{B_n(r)}=0$.
\enddemo

\subhead 6.\ \ Bootstrapping\endsubhead

We continue using $\|\quad \|$ for $l^1$--norm on $\bC^N$.
Our first application of Lemma 5.1 gives a slight improvement on Corollary 3.2.

\proclaim{Proposition 6.1} Given $0 < r < R$, there is a constant $C > 1$ with the following property.
Suppose that $f$ is a complex valued, bounded and measurable $(0,1)$
form on the ball $B_N(R)\subset \bC^N$.
If $\opartial f=0$ then there is a $u\in C(B_N(r))$ such that $\opartial u=f$ and for $z\in B_N(r)$
$$
|u(z)|\leq R C^{\# z} |f|_0.\tag6.1
$$
$C$ depends only on the ratio $r/R < 1$ but not on $N$.
\endproclaim

\demo{Proof} Put $r'=(R+r)/2$, and let $u\in C(B(r'))$ denote the canonical solution constructed in Lemma 5.1, Corollary 5.2 based on the nodes
$Z(0)=0,\ Z(k)=r' k/\|k\|,\quad k\not= 0$.
We claim that (6.1) holds.
Indeed, suppose $\#z=n > 0$; we can assume that the first $n$ coordinates of $z$ are nonzero.
By virtue of Corollary 3.2 there is a continuous function $U$ on $B_n(R)=\bC^n\cap B_N(R)$ that solves $\opartial U=F|_{B_n(R)}$ and satisfies
$$
|U(\zeta)|\leq 2(1+2\sqrt{n}) R \left({R\over R-\|\zeta\|}\right)^n |f|_0
\leq R\left({6R\over R-\|\zeta\|}\right)^n |f|_0,
$$
$\zeta\in B_n(R)$.
Thus the hypothesis of Lemma 5.1 is satisfied with 
$Q=6R/(R-r')$, $A=R|f|_0$ and (5.5) gives
$$
|u(z)|\leq |U(z)|+R|f|_0 \Delta (Q, z/r')\leq R(Q^{\# z}+e^{c\# z} Q^{\# z}) |f|_0,
$$
by Corollary 4.3 if $z\in B_N(r),\ z\not= 0$.
On the other hand (5.3) implies $u_k(0)=0$ if $k\not= 0$; also $u_0(0)=0$ by
the normalization condition $u_k(Z(k))=0$; therefore $u(0)=\Sigma u_k (0)=0$.
Hence (6.1) follows.

In order to take more advantage of Lemma 5.1 and Corollary 5.2, especially of the inequality (5.5), we shall need a slightly stronger statement than Proposition 6.1:\ in addition to being able to solve $\opartial u=f$ with a $u=U$ that satisfies an estimate like (6.1) we want to have $U(Z)$ under control for an arbitrary but fixed $Z\in B_N(r)$.
The next proposition does just that:

\proclaim{Proposition 6.2} Given $0 < r < R$ there is a constant $Q > 1$ depending only on $r/R < 1$ with the following property.
Suppose $Z\in B_N(r)$, and $f$ is a complex valued, Lipschitz continuous $(0,1)$ form on $B_N(R)$.
If $\opartial f=0$ then there is a function $U\in C^1 (B_N (r))$ such that $\opartial U=f$, and for $z\in B_N(r)$
$$
|U(z)|\leq \|z-Z\| Q^{1+\#z} (|f|_0+R|f|_1).\tag6.2
$$
\endproclaim

\demo{Proof} The Proposition is clearly true when $N=1$; we shall prove it for all $N$ by induction.
Thus, assume it is true with $N$ replaced by $N-1$, and we shall also assume without loss of generality that $Z_\nu\geq 0,\ \nu=1,\ldots N$.
\enddemo

With $R'=R-Z_N$ define a holomorphic mapping $\pi\colon B_N(R)\to B_{N-1} (R')$ by
$$
\pi (z)=z' (R-Z_N)/(R-z_N)\quad ,\ z=(z',z_N).
$$
If $\var\colon \bC^{N-1}\to \bC^N$ denotes the embedding $\var(z')=(z',Z_N)$, then $\var\circ\pi$ is a central projection with center $(0,\ldots,0,R)$.
Also $\pi\circ \var=$ id.
One checks that $\pi(B_N(r))\subset B_{N-1}(r'),\ r'=(R-Z_N) r/R$.
In particular $r'/R'=r/R$ so we can apply our inductive hypothesis with
the form $f'=\var^* f$ to obtain $U'\in C^1(B_{N-1}(r'))$ satisfying
$\opartial U'=f'$ and
$$
|U'(z')|\leq \|z'-Z'\|Q^{1+\#z'} (|f'|_0+R'|f'|_1)\quad,\quad Z=(Z',Z_N).\tag6.3
$$
Define a closed $(0,1)$ form $F=\sum_1^N F_\nu d\overline z_\nu=f-\pi^* f'$
on $B_N(R)$.
Then $F=f-\opartial(\pi^* U')$ on $B_N(r)$ and $\var^* F=0$ so that
$$
F_\nu (z', Z_N)=0\quad ,\ \nu=1,\ldots,N-1.\tag6.4
$$
Also
$$
\align
F&=\textstyle{\sum^{N-1}_1} F_\nu d\oz_\nu-(\oz_N-\oZ_N)\opartial F_N+\opartial ((\oz_N-\oZ_N) F_N)\\
&=(z_N-Z_N) g+\opartial((\oz_N-\oZ_N) F_N),
\endalign
$$
where
$$
g={\sum_1^{N-1} F_\nu d\oz_\nu\over z_N-Z_N}-{\oz_N-\oZ_N\over z_N-Z_N}\opartial F_N\tag6.5
$$
if $z_N\not= Z_N,\ g=0$ if $z_N=Z_N$.

Put $R_0=(R+r)/2$.
Lipschitz continuity of $f$, therefore of $F$, and (6.4) imply that $g$ is a bounded, measurable, closed $(0,1)$ form on $B_N(R_0)$, hence by Proposition 6.1 a $v\in C(B_N(r))$ can be found that solves
$\opartial v=g$ and satisfies
$$
|v(z)|\leq R C^{\#z}|g|_{0,\ B_N(R_0)}\tag6.6
$$
with some $C > 1$ depending only on $r/R$.
Therefore
$$
U=\pi^* U'+(z_N-Z_N)v+(\overline z_N-\overline Z_N) F_N
$$
solves $\opartial U=f$ on $B_N(r)$.
By elliptic regularity, cf.~Proposition 2.3, $U$ is $C^1$.
To complete the proof we must estimate $U(z)$.

First observe that when $z\in B_N(r)$
$$
\align
\|\pi(z)-Z'\|&=\Bigg\| {R-Z_N\over R-z_N} z'-Z'\Bigg\|\leq\\
&\leq \|z'-Z'\|+\bigg|{z_N-Z_N\over R-z_N}\bigg| \|z'\|\leq \|z'-Z'\|+{r\over R}|z_N-Z_N|,
\endalign
$$
so that (6.3) gives
$$
|U' (\pi(z))|\leq (\|z'-Z'\|+{r\over R}|z_N-Z_N|) Q^{1+\#z} (|f|_0+R|f|_1).\tag6.7
$$

Next we estimate $|F|_0$ and $|F|_1$.
One first computes
$$
|\pi|_{1,B_N(R_0)}\leq 2R^2/(R-R_0)^2\quad,\quad
|d\pi|_{1, B_N(R_0)}\leq 4R^3/(R-R_0)^4;
$$
next, by (2.4)
$$
\align
|\pi^* f'|_{0,B_N(R_0)}&\leq {2R^2\over (R-R_0)^2} |f|_0,\\
|\pi^* f'|_{1,B_N(R_0)}&\leq {4R^3\over (R-R_0)^4} |f|_0+{4R^4\over (R-R_0)^4} |f|_1.
\endalign
$$
It follows that $|F|_{0,B_N(R_0)}$ resp.~$|F|_{1,B_N(R_0)}$ are dominated by twice the right hand sides above.
(6.4) and (6.5) imply
$$
|g|_{0,B_N(R_0)}\leq {256 R^3\over (R-r)^4} |f|_0+
{256R^4\over (R-r)^4} |f|_1,\text{ and}
$$
$$
|v(z)|\leq 256 C^{\#z}\left\{{R^4\over (R-r)^4}|f|_0+
{R^5\over (R-r)^4} |f|_1\right\}
$$
by (6.6).
Since $|F_N|_{0, B_N(r)}\leq |F|_{0,B_N(R_0)}\leq 16R^2|f|_0/(R-r)^2$,
combining these last two inequalities with (6.7) we obtain, when $z\in B_N(r)$
$$
\align
|U(z)|&\leq \left(\|z'-Z'\|+{r\over R}|z_N-Z_N|\right)Q^{1+\#z} (|f|_0+R |f|_1)\\
&+256 |z_N-Z_N| C^{\#z}\left\{{R^4\over (R-r)^4}|f|_0+{R^5\over (R-r)^4}|f|_1\right\}\\
&+16 |z_N-Z_N|{R^2\over (R-r)^2} |f|_0\\
&\leq \|z-Z\| Q^{1+\#z} (|f|_0+R|f|_1),
\endalign
$$
provided $Q\geq 272 CR^5/(R-r)^5$, q.e.d.
\enddemo

\subhead 7.\ \ The equation $\opartial u=f$ in a ball in $l^1$\endsubhead

We are ready to tackle the $\opartial$ equation in $l^1$.
In this section we prove

\proclaim{Theorem 7.1} Given $0 < r < R$ there is a constant $Q > 1$, depending only on the ratio $r/R$, with the following property.
If a Lipschitz continuous $(0,1)$ form $f$ on $B(R)$ is $\opartial$--closed then the equation $\opartial u=f$ has a solution $u\in C^1 (B(r))$.
This solution can be chosen so that
$$
|u(z)|\leq 2RQ\Delta (Q, z/r)(|f|_0+R|f|_1),\ z\in B(r).\tag7.1
$$
\endproclaim

\demo{Proof} Embed the spaces $\bC^N$ into $l^1$ by the map $(z_\nu)^N_1\mapsto (z_1,\ldots,z_N,0,\ldots)$, and let $\pi_N\colon l^1\to \bC^N\subset l^1$ denote the projection $(z_\nu)_1^\infty\mapsto (z_1,\ldots,z_N,0,\ldots)$.
We shall first prove the theorem assuming that $f=\pi_N^* g$ with some $g\in C_{0,1} (B_N(R))$.
\enddemo

Let $v\in C(B_N (R))$ denote the canonical solution of the equation
$\opartial v=g$ constructed in Lemma 5.1, Corollary 5.2, based on the nodes $Z(k)=rk/\|k\|$.
We can estimate $v(Z)$ for arbitrary $Z\in B_N(r)$ as follows.
Put $r'=(R+r)/2$, and apply Proposition 6.2 with $r$ replaced by $r'$.
We obtain a $Q > 1$ depending only on $r/R$, and a $U\in C^1(B_N (r'))$ that solves $\opartial U=g$ and satisfies
$$
|U(z)|\leq \|z-Z\| Q^{1+\#z}(|g|_0+R|g|_1),\ z\in B_N (r').
$$
In particular (5.4) in Lemma 5.1 is satisfied with $A=2RQ(|g|_0+R|g|_1)$.
Therefore the lemma gives the following estimate for the canonical solution $v$:
$$
|v(Z)|=|v(Z)-U(Z)|\leq 2RQ\Delta (Q, Z/r)(|g|_0+R|g|_1).
$$
Thus $u=\pi_N^* v$ is as required.

Now for a general $f$ as in the theorem define $f^N=\pi_N^* f=\pi_N^* (f|_{\bC^N})$.
Since $\pi_N\colon l^1\to l^1$ converges to the identity as $N\to\infty$, uniformly on compacts, $f^N\to f$ uniformly on compact subsets of $B(R)\times l^1$.
Let $u^N\in C(B(r))$ denote the solution of the equation $\opartial u^N=f^N$ 
constructed in the first half of the proof. We will show the sequence $u^N$
is convergent on $B(r)$. Corollary 5.3 implies 
$u^{N+1}|_{B_N(r)}=u^N|_{B_N(r)}$, hence the $u^N$ converge on the dense
subset $\{z\in B(r) : \#z<\infty\}$. Furthermore, (7.1) is satisfied,
with $N$'s appended; in particular, the $u^N$ are locally uniformly bounded,
hence locally equicontinuous by Lemma 4.1 and Proposition 2.4.
It follows the sequence $u^N$ is uniformly convergent on compact subsets of $B(r)$.
The limit $u$ satisfies (7.1); it is continuous by Proposition 2.1; by Proposition 2.2 $\opartial u=f$, and by Proposition 2.3 $u\in C^1(B(r))$.
With this the proof of Theorem 7.1 is complete.

\subhead 8.\ \ The proof of Theorem 1.1\endsubhead

Theorem 1.1 will be derived from Theorem 7.1 and an approximation theorem.
This latter concerns approximating holomorphic functions on $B(R)$ by entire functions; the approximation should take place on $B(r)$ for some $r < R$.
Ideally, one would approximate uniformly on $B(r)$.
However, as indicated in the Introduction, such a strong approximation theorem is not known, and we shall use a norm different from sup norm on $B(r)$ to quantify approximation.

Thus, let $\Sigma a_k z^k$ be the monomial expansion of a function $h$ holomorphic on some neighborhood of $0\in l^1$ (see Section 4), and with $0 < r < \infty$ put
$$
[h]_r=\sup_k |a_k| r^{\|k\|} k^k/\|k\|^{\|k\|}\leq \infty.
$$
Note that $r'\geq r$ implies $[h]_{r'}\geq [h]_r$.

\proclaim{Theorem 8.1} Suppose $\varphi$ is a holomorphic function on $B(R)\subset l^1$.
For any $0 < r < R$ and $\var > 0$ there is a function $\psi$, holomorphic on $l^1$, such that $[\varphi-\psi]_r\leq \var$.
\endproclaim

The theorem is equivalent to the following 

\proclaim{Lemma 8.2}
Let $0 < r < R$, and let $\calK$ be a sequence of integer multiindices $k=(k_1,k_2,\ldots) \geq 0,\ \#k < \infty$.
Suppose for each $k\in \calK$ we are given $c_k\in \bC$ such that
$$
\inf_\calK |c_k| r^{\|k\|} > 0;\text{ and}\tag8.1
$$
$$
\lim_\calK c_k \tau^k=0\tag8.2
$$
whenever $\tau=(\tau_\nu)$ satisfies
$$
0 \leq \tau_\nu < R\quad ,\qquad \lim_\nu\tau_\nu=0.\tag8.3
$$
Then (8.2) holds for any $\tau=(\tau_\nu)$ satisfying $\lim\limits_\nu \tau_\nu=0$.
\endproclaim

\demo{Proof} Supposing $\lim \tau_\nu=0$, with $0 < \delta < 1$ put $\sigma_\nu=\tau_\nu^\delta r^{1-\delta}\to 0\quad (\nu\to\infty)$.
Fix $\delta$ so that $\sigma_\nu < R$ for all $\nu$.
This can be arranged, since when $\nu$ is sufficiently large $\tau_\nu < r$, so that $\sigma_\nu < R$ however $\delta$ is chosen; and letting $\delta$ be sufficiently small will take care of the remaining finitely many $\nu$'s.
With such $\delta$ and $\sigma=(\sigma_\nu)$ we have $\lim_\calK c_k \sigma^k=0$ by (8.2), whence $\tau_\nu=\sigma_\nu^{1/\delta} r^{1-1/\delta}$ gives
$$
c_k \tau^k=(c_k \sigma^k)^{1/\delta} (c_k r^{\|k\|})^{1-1/\delta}\to 0.
$$ 
\enddemo

\demo{Proof of Theorem 8.1} If $\varphi(z)=\Sigma_k \alpha_k z^k$ is the monomial expansion of $\varphi$, let
$$
\calK=\{k\colon |\alpha_k| r^{\|k\|} k^k/\|k\|^{\|k\|}\geq \var\}\quad ,\quad \text{ and }\quad\psi(z)=\sum\limits_{k\in\calK} \alpha_k z^k.
$$
Lemma 8.2 and Theorem 4.5 together imply $\psi$ is holomorphic on $l^1$; clearly also $[\varphi-\psi]_r\leq \var$, as claimed.
\enddemo

\demo{Proof of Theorem 1.1} Choose a sequence $0 < r_1 < r_2 < \ldots \to R$, and for each $n=1,2,\ldots$, let $u_n\in C^1(B(r_n))$ solve $\opartial u_n=f|_{B(r_n)}$, cf.~Theorem 7.1.
Thus $u_n-u_{n-1}$ is holomorphic on $B(r_{n-1})$; by Theorem 8.1 there are entire functions $\psi_n$ such that $[u_n-u_{n-1}-\psi_n]_{r_{n-2}}\leq 2^{-n},\ n=2,3,\ldots$.
With $v_n=u_n-\Sigma_2^n \psi_j$ it follows that for $0 < r < R$
$$
\lim_{n,p\to \infty} [v_n-v_p]_r=0,
$$
and so $v_n$ uniformly converges on compact subsets of $B(r)$ by Theorem 4.4(b).
This being true for any $r < R$, in view of Proposition 2.1 $v_n$ converges to some $u\in C(B(R))$, uniformly on compacts in $B(R)$.
Finally $f=\lim\opartial u_n=\lim\opartial v_n=\opartial u$ by Proposition 2.2; thus $u\in C^1 (B(R))$, resp.~$u\in C^m (B(R))$ if $f\in C^m_{0,1} (B(R))$, by Proposition 2.3.
\enddemo

\subhead 9.\ \ An example\endsubhead

Although in Theorem 1.1 one can very likely relax the hypothesis on the smoothness of $f$, in $l^1$ for the solvability of the equation $\opartial u=f$ it is not enough to assume that $|f|_0$ alone is finite --- in contrast with the finite dimensional situation.
Below we will give examples in $l^p$ spaces of closed $(0,1)$ forms $f$ of higher and higher regularity that are not exact.
In this section $\|\quad\|$ will denote the norm on an $l^p$ space and $B(r)=\{z\in l^p\colon \|z\| < r\}$.

\proclaim{Theorem 9.1} For any $p=1,2,\ldots$ there is a $\opartial$--closed $f\in C_{0,1}^{p-1}(l^p)$ such that on no open set $\Omega\not= \emptyset$ is the equation $\opartial u=f$ solvable.
Moreover, $f$ can be chosen to satisfy the growth condition
$$
|f|_{0,B(r)}\leq r^{p-1},\quad r > 0.\tag9.1
$$
\endproclaim

The construction to be outlined below is but a minor extension of one given by Coeur\'e for the case $p=2$, see [C,M].
Aside from its contrast with Theorem 1.1, Theorem 9.1 is of interest because we have proved elsewhere that a growth condition like (9.1) and the assumption $f\in C_{0,1}^p (V)$ --- rather than $C^{p-1}$ --- do imply the solvability of the equation $\opartial u=f$ on any (sequentially complete) locally convex topological vector space $V$, see [L, Theorem 9.1] for a precise formulation.

\proclaim{Proposition 9.2} Suppose a compactly supported function $\varphi\in C^{p-1}(\bC)$ satisfies $|\varphi(\zeta)|\leq |\zeta|^{p-1},\ \zeta\in\bC$.
Then
$$
f(z;\xi)=\sum_{\nu=1}^\infty\varphi(z_\nu)\overline{\xi_\nu}\tag9.2
$$
defines a $\opartial$-closed $f\in C_{0,1}^{p-1}(l^p)$, and
$$
|f|_{0,B(r)}\leq r^{p-1}.\tag9.3
$$
\endproclaim

\demo{Proof} By H\"older's inequality (9.2) converges for $z,\xi\in l^p$, and $|f(z;\xi)|\leq \Sigma |z_\nu|^{p-1} |\xi_\nu|\mathbreak
\leq \|z\|^{p-1} \|\xi\|$,
so that (9.3) is satisfied.
Similarly, the series gotten by taking iterated derivatives in (9.2) in the directions $\eta^1,\ldots,\eta^{p-1}\in l^p$ is dominated by the series
$$
\text{const }\sum_{\nu=1}^\infty|\xi_\nu|\prod^{p-1}_{j=1}|
\eta_\nu^j|\leq\text{ const } \|\xi\| \prod\limits_j\|\eta^j\|
$$
(this last inequality follows by $p$ applications of H\"older's inequality), hence it is convergent and its sum is locally uniformly bounded on $l^p\times l^p\times \ldots l^p$.
Propositions 4.2 and 2.1 imply that this dominating series is uniformly 
convergent on compacts in $l^p\times l^p\times \ldots\times l^p$ and
therefore $f\in C^{p-1}_{0,1}(l^p)$.
Since the partial sums of (9.2) represent the forms $\sum_1^n \varphi(z_\nu)d\overline z_\nu$, which are manifestly $\opartial$--closed, their limit $f$ is also $\opartial$--closed.

\proclaim{Proposition 9.3} The function
$$
{\partial\over \partial\ozeta} \zeta^p \log \log |\zeta|^{-2}={\zeta^p\over \ozeta\log |\zeta|^2}
$$
is of class $C^{p-1}$ on the unit disc in $\bC$.
\endproclaim

We omit the proof, a simple exercise.
The point here is that $\zeta^p\ \log \log |\zeta|^{-2}$ is not $C^p$ while its $\opartial$ is $C^{p-1}$, a subtle phenomenon well understood in harmonic analysis.
This mild irregularity of the one dimensional $\opartial$ operator will have as consequence nonsolvability in infinite dimensions:

\demo{Proof of Theorem 9.1} Let $\lambda$ be a compactly supported function on $\bC$ that is $C^p$ off 0 and agrees with $\zeta^p\log\log |\zeta|^{-2}$ in a neighborhood of 0.
Arrange also that $\varphi(\zeta)=\partial\lambda(\zeta)/\partial\ozeta$ is dominated by $|\zeta|^{p-1}$, and construct $f$ as in Proposition 9.2.
In view of Propositions 9.2, 9.3 all we have to do to prove the theorem is to show that on no open set is there a continuous function $u$ with $\opartial u=f$.
\enddemo

Assume first that on a ball $B(2R)\subset l^p$ there is a solution $u$.
Upon shrinking $R$ we can assume that $u$ is bounded on $B(2 R)$ by a number $M$ and also that $\lambda(\zeta)=\zeta^p\log\log |\zeta|^{-2}$ when $|\zeta| < 2R$.
For arbitrary $N=1,2,\ldots$ embed $\bC^N$ into $l^p$ by $\bC^N\ni w\mapsto (w,0,0\ldots)\in l^p$.
Since $\opartial u|_{\bC^N}=f|_{\bC^N}=\opartial \sum^N_{\nu=1} \
\lambda(z_\nu)$, it follows that $u(z)-\sum^N_1\lambda (z_\nu)$ is holomorphic on $B_N(2R)$.
Thus we have found a holomorphic function $h$ such that
$$
\bigg|\sum^N_1\lambda(z_\nu)-h(z_1,\ldots,z_N)\bigg|\leq M,\ z\in B_N(2R).\tag9.4
$$

Next we claim that there is even a homogeneous polynomial $h$ of degree $p$ that satisfies (9.4).
Indeed, all we have to do is to replace the original $h$ by
$h'(z)=\int_0^1 e^{-2\pi pit}h(e^{2\pi it} z)\ dt$ and note that
$$
\int_0^1 e^{-2\pi pit}\sum\limits_\nu \lambda (e^{2\pi it} z_\nu)\ dt=
\sum\limits_\nu \lambda_\nu (z_\nu).
$$
In the same spirit, with $\rho$ the action of the torus $T^N$ on $\bC^N$ as before let $G\subset T^N$ denote the subgroup of elements of order $p,\ G\cong (\bZ/p)^N$.
Restrict $\rho$ to $G$; upon averaging (9.4) over orbits of $\rho|_G$ we obtain a polynomial $h$ satisfying (9.4) that is invariant under $\rho|_G$, i.e., $h$ is a linear combination of the monomials $z_\nu^p$.
Finally, symmetrization in the variables $z_\nu$ yields an $h$ of form $a\sum_1^N z_\nu^p$ that satisfies (9.4), with $a=a_N$ a constant.

With an arbitrary $n\leq N$ put
$$
z_\nu=\cases Rn^{-1/p}&\text{$\nu=1,\ldots,n$}\\
0&\text{$\nu=n+1,\ldots, N$}\endcases.
$$
Thus $z=(z_\nu)\in B_N(2R)$.
Substitute this $z$ into (9.4) to obtain
$$
|R^p\log\log (R^{-2} n^{2/p})-a_N R^p|\leq M,
$$
which cannot hold for all $n\leq N$ if $N$ is sufficiently large.
This contradiction shows that $\opartial u=f$ is not solvable on any neighborhood of the origin.

Now we conclude as follows.
Suppose $\Omega\subset l^p$ is a nonempty open set, and select a point $Z\in\Omega$ such that $Z_\nu=0$ for $\nu > N$.
Embed $l^1$ into $l^1$ by the map $z\mapsto E(z)=w$, where $w_\nu=Z_\nu$ if $\nu\leq N$, $w_\nu=z_{\nu-N}$ if $\nu > N$.
If $f$ were exact in a neighborhood of $Z=E(0)$ then $E^*f=f$ would also be exact in a neighborhood of 0.
Since this is not the case, $\opartial u=f$ is not solvable on $\Omega$.
\documentstyle{amsppt}

\bigskip
\baselineskip=12pt
\noindent
L\'aszl\'o Lempert\newline
Department of Mathematics\newline
Purdue University\newline
West Lafayette, IN\ \ 47907--1395

\bye